\documentclass[12pt]{article}
\usepackage[utf8]{inputenc}
\usepackage{authblk}
\title{Extended Set Difference : Inverse Operation of Minkowski Summation}
\author{
  Arie Beresteanu\footnote{University of Pittsburgh, Pittsburgh, PA, USA. e-mail: arie@pitt.edu.}
  \qquad
  Behrooz Moosavi Ramezanzadeh\footnote{University of Pittsburgh, Pittsburgh, PA, USA. e-mail: behroozmoosavi@pitt.edu.}
}
\date{November 2024}
\usepackage{pgfplots} 
\usepackage{caption}
\usepackage{amssymb}
\usepackage{amsfonts}
\usepackage{amsmath}
\usepackage{scalefnt}
\usepackage{fancybox}
\usepackage[ruled,vlined]{algorithm2e}  % Use algorithm2e package
\usepackage{float}
\usepackage{natbib}
\usepackage{theorem}
\usepackage{lscape}
\usepackage{hyperref}
\usepackage{algpseudocode}
\usepackage{graphicx}
\usepackage{appendix}
\usepackage{caption}
\usepackage{subcaption}
\usepackage{hyperref}

\usepackage{cite}
\usepackage{tikz}
\usetikzlibrary{shapes.geometric}
\usepackage{comment}
\newtheorem{theorem}{Theorem}
\newenvironment{theorem*}{\vskip 0.1in \noindent {\bf Theorem \/}}{\vskip 0.1in}

\newtheorem{definition}{Definition}
\newtheorem{example}{Example}

\newtheorem{lemma}{Lemma}

\newtheorem{proposition}[theorem]{Proposition}
{\theorembodyfont{\upshape} \newtheorem{remark}{Remark}}

\newenvironment{proof}[1][Proof]{\textbf{#1.} }{\ \rule{0.5em}{0.5em}}
\renewcommand{\cite}{\citet}
\usepackage{pgfplots}
\usepackage{pgfmath}
\pgfplotsset{compat=1.17}
\newcommand{\Rd}{\mathbb{R}^d}
\newcommand{\Sd}{\mathbb{S}^{d-1}}
\newcommand{\rin}{r_{in}}
\newcommand{\Rout}{R_{out}}
\newcommand{\KCn}{\mathcal{K}^d_{kc}}
\definecolor{linkcolor}{rgb}{0,0,1} % Blue color for links.
\hypersetup{
    colorlinks=true,
    linkcolor=linkcolor,
    filecolor=magenta,      
    urlcolor=linkcolor,
}

\usepackage{pgfplots}
\usepackage{pgfmath}
\pgfplotsset{compat=1.17}
\date{\today}
\begin{document}
\maketitle

\begin{abstract}
This paper introduces the extended set difference, a generalization of the Hukuhara and generalized Hukuhara differences, defined for compact convex sets in $\mathbb{R}^d$. The proposed difference guarantees the existence for any pair of such sets, offering a broader framework for set arithmetic. The difference is not necessarily unique, but we offer a bound on the variety of solutions. We then show how a unique solution is found as the limit of strictly convex perturbations of the concept of set difference. The definition of the extended set difference is formulated through an optimization problem, which provides a constructive approach to its computation. The paper explores the properties of this new difference, including its stability under orthogonal transformations and its robustness to perturbations of the input sets. We propose a computational algorithm based on linear optimization to efficiently calculate the extended set difference.
\end{abstract}

\noindent\textbf{Credit authorship contribution statement: }

\textbf{Arie Beresteanu: } Writing - review \& editing, Writing - original draft, Methodology, Investigation, Conceptualization. 

\textbf{Behrooz Moosavi Ramezanzadeh: } Writing - review \& editing, Writing - original draft, Methodology, Investigation, Conceptualization.

\noindent\textbf{Conflict of Interest:} The authors declare that they have no conflict of interest.

\break

\break

%%%%%%%%%%%%%%%%%%%%%%%%%%%%%
%%%%%   Introduction    %%%%%
%%%%%%%%%%%%%%%%%%%%%%%%%%%%%
\section{Introduction}

The Minkowski summation of two subsets of \(\mathbb{R}^d\), defined as \(A \oplus B = \{a+b : a \in A, b \in B\}\), does not have an inverse operation. For three vectors \(a, b, x \in \mathbb{R}^d\), equation \(a - b = x\) is equivalent to \(a = b + x\). Therefore, \cite{hukuhara-1967} suggested defining the difference \(A \ominus_H B\) as a set \(X\) such that \(A = B \oplus X\). However, \(A \ominus_H B\) does not necessarily exist. A necessary condition for the existence of \(A \ominus_H B\) is that \(\exists x \in \mathbb{R}^d\) is such that \(x + B \subset A\) where $x + B$ means $\{x\} \oplus B$. \cite{STEFANINI2010} generalizes the Hukuhara difference by observing that in \(\mathbb{R}^d\), \(a - b = x\) is equivalent to \(a + (-x) = b\). Thus, \cite{STEFANINI2010} defines \(A \ominus_g B\) to be equal to \(A \ominus_H B\) if the regular Hukuhara difference exists, and to a set \(X\) such that \(A \oplus (-X) = B\) if the regular Hukuhara difference does not exist.\footnote{\(-X = \{-x : x \in X\}\) is the opposite set of \(X\).} A necessary condition for the generalized difference to exist is that \(\exists x \in \mathbb{R}^d\) such that either \(x + B \subset A\) or \(x + A \subset B\). For example, if \(A = [0,1]\) and \(B = [2,4]\), then \(\nexists X\) such that \([2,4] \oplus X = [0,1]\). On the other hand, \([0,1] \ominus_g [2,4] = [-3,-2]\). For intervals in $\mathbb{R}^1$ the generalized Hukuhara difference always exists and is unique. 

The translation condition in the generalized Hukuhara difference is sufficient to guarantee the existence of the set difference for intervals in $\mathbb{R}$. However, this condition is not sufficient in $\Rd$ for $d \geq 2$. For example, if \(A\) is the unit cube and \(B\) is the unit circle in \(\mathbb{R}^d\) for \(d \geq 2\), we cannot find a set \(X \subset \mathbb{R}^d\) such that \(A = B \oplus X\), even though \(B \subset A\).

The purpose of this paper is to introduce an alternative set difference that extends the generalized Hukuhara difference and is well defined for any two compact convex sets in \(\mathbb{R}^d\) and any \(d \in \mathbb{N}\). This new difference is formulated through an optimization problem that minimizes a criterion function based on the Hausdorff distance between $A$ and $B\oplus X$, rather than directly trying to achieve $A = B \oplus X$. The uniqueness of the minimizer, $X$ is not necessarily guaranteed. We provide an upper bound on the variety of solutions and show that two potential solutions cannot be far away from each other in the Hausdorff sense. We then show that by a refinement of the difference concept through a series of strictly convex perturbations of the criteria function, a unique solution can be achieved. Despite the non uniqueness, an approximation to a set $X$ which minimizes the Hausdorff distance between $A$ and $B\oplus X$ can be computed using a Linear Programming problem. We provide an algorithm to implement this approximation.

%%%%%   Notations   %%%%%
\subsection{Notations}
We denote by $\KCn$ the set of all compact and convex sets in the Euclidean space $\left<\mathbb{R}^d,||\cdot||_2\right>$. For convenience, we omit the subscript and use $||\cdot||$ for the Euclidean norm. We let $\Sd = \{u \in \Rd : ||u|| = 1 \}$ denote the unit circle in $\Rd$. For a set $A \in \KCn$ we denote by $h_A:\Sd \rightarrow \mathbb{R}$ its support function $h_A(u)=sup\{u\cdot a\in A\}$(see \cite{rockafellar-2015}). For two sets $A,B \in \KCn$, $d_H(A,B)$ denotes the Hausdorff distance defined as $d_H(A, B) = \max \left\{ \sup_{a \in A} \inf_{b \in B} d(a, b), \sup_{b \in B} \inf_{a \in A} d(b, a) \right\}$. Minkowski summation of two sets $A,B \in \KCn$ is denoted $A \oplus B$ and defined as $A \oplus B =\{a+B : a\in A, \ b \in B \}$. The set $\mathbb{B}$ denotes the unit ball in $\Rd$ centered at the origin. For $r\in \mathbb{R}_+$, $\mathbb{B}_r= r\mathbb{B}$ is a ball with radius $r$ centered at the origin and $\mathbb{B_r}(c) = r\mathbb{B}\oplus{c}$ is a ball of radius $r$ centered at a point $c\in \Rd$.

\par
\subsection{The structure of the paper}
In Section \ref{sec:ExistUnique}, we define the extended set difference as a solution to an optimization problem and show that it exists for any two compact convex sets in $\KCn$. In Section \ref{sec:Properties} we provide basic properties of the extended set difference and investigate its behavior. We also provide a bound on the distance between any two solutions to this optimization problem. A unique solution through a strictly convex pertubation is offered in Section \ref{sec: refinement}. In Section \ref{sec:computation} we show that the extended set difference we propose can be approximated arbitrarily closely by a convex polygon and be computed using an LP problem. 

%%%%%%%%%%%%%%%%%%%%%%%%%%%%%%%%%%%%%%
%%%%%  Extended Set Difference   %%%%%
%%%%%%%%%%%%%%%%%%%%%%%%%%%%%%%%%%%%%%
\section{ Extended Set Difference}\label{sec:ExistUnique}
%%%%%   Existance    %%%%%
\subsection{Existence}

While $A \ominus_g B$ exists for all \(A, B \in \mathcal{K}^1_{kc}\), it does not generally exist for \(A, B \in \KCn\) when \(d \geq 2\). To address this lack of existence, we extend the framework of set-valued arithmetic by generalizing the Hukuhara difference (\cite{hukuhara-1967}) and the generalized Hukuhara difference (\cite{STEFANINI2010}). We propose the following alternative definition for set difference.

\begin{definition}[Extended Set Difference]\label{def:extended-diff}
For \( A, B \in \KCn \), the extended set difference is the collection
\begin{equation}
A \ominus_e B = \arg \inf_{X \in \KCn} d_H(A, B \oplus X),
\end{equation}
where \( B \oplus X = \{ b + x : b \in B, x \in X \} \) is the Minkowski sum, and \( d_H(A, B) = \max \{ \sup_{a \in A} \inf_{b \in B} \|a - b\|, \sup_{b \in B} \inf_{a \in A} \|a - b\| \} \) is the Hausdorff distance. The minimal distance achieved is denoted \( \delta = \inf_{X \in \KCn} d_H(A, B \oplus X) \).
\end{definition}

We seek a set $X$ that minimizes the Hausdorff distance between \(A\) and \(B \oplus X\) rather than find a set $X$ such that $A=B\oplus X$. Generally, there can be more than one set in $\KCn$ that minimizes the criteria function in (\ref{def:extended-diff}). In other words, $A \ominus_e B \subset \KCn$. We therefore call $A \ominus_e B$ the extended set difference collection. Theorem \ref{thm:Existance} below ensures that \(A \ominus_e B\) is always well defined in the sense that $A \ominus_e B$ is not an empty collection. Since (\ref{def:extended-diff}) may have more than one solution, Theorem \ref{thm:Scope} in section \ref{sec:Properties} establishes a bound on the variety of solutions.

\begin{theorem}[Existence]\label{thm:Existance}
Let \(A, B \in \KCn\). Then there exists a set \(X \in \KCn\) such that $X \in A \ominus_e B.$
\end{theorem}
\begin{proof}
Given $A,B \in \KCn$, let $f: \KCn \rightarrow \mathbb{R}$ be \(f(X) = d_H(A, B \oplus X)\). We need to show that $f()$ attains its minimum in \(\KCn\).

Since $\KCn$ is close under Minkowski addition, $B, X \in \KCn$, implies $B \oplus X\ \in \KCn$. Moreover, the function \(f(X)\) is continuous with respect to \(X\). 
Specifically, if $\{X_n\}$ is a sequence of sets in $\KCn$ and $X \in \KCn$ such that $d_H(X_n,X) \rightarrow 0$ as $n \rightarrow \infty$, then, by Theorem 1.1.2 in \cite{li-2002}, $\lim_{n\to\infty} f(X_n) = f(X).$

For any set $E \in \KCn$, define $\KCn|E = \{K \cap E : K \in \KCn \}$ to be the collection of compact convex sets in $\Rd$ contained in the set $E$. Since $A,B$ are compact, there is a ball $\mathbb{B}_r$ with a finite $r>0$ such that $A \subset B \oplus \mathbb{B}_r$. Therefore, the infimum in (\ref{def:extended-diff}) can be taken over $\KCn|\mathbb{B}_r$ rather than $\KCn$. blaschke Selection Theorem (see \cite{willmore-1952}), states that every bounded sequence of compact convex sets in $\KCn$ has a convergent subsequence whose limit is also in $\KCn$. Moreover, this limit is bounded in $\mathbb{B}_r$ as well. Therefore, $\KCn|\mathbb{B}_r$ is sequentially compact. 
$\left(\KCn|\mathbb{B}_r,d_H \right)$ is a metric space and therefore sequential compactness is equivalent to compactness by the Heine-Borel Theorem.

Thus, there exists \(X^{\ast} \in \KCn\) such that
\[
f(X^{\ast}) = \inf_{X \in \KCn} f(X) = \inf_{X \in \KCn} d_H(A, B \oplus X).
\]
Moreover, \(X^{\ast} \in A \ominus_e B  \) is, by construction, compact and convex, and the Minkowski sum \(B \oplus X^{\ast}\) remains compact. This shows that \(A \ominus_e B\) is not empty.
\end{proof}

In what follows we denote by $\mathcal{M}$ the set of solutions for the minimization problem in definition \ref{def:extended-diff}. We do so without referring to the sets $A$ and $B$. Theorem \ref{thm:Existance} shows that $\mathcal{M}\neq \emptyset$. In the next section we explore the properties of the set difference operator and of the set of solutions $\mathcal{M}$.

%%%%%%%%%%%%%%%%%%%%%%%%%%%%%%%%%%%%%%
%%%% Properties of set difference %%%%
%%%%%%%%%%%%%%%%%%%%%%%%%%%%%%%%%%%%%%

\section{Properties of Set Differences\label{sec:Properties}}
The collection of sets that solve the optimization problem in (\ref{def:extended-diff}) is non-empty by Theorem \ref{thm:Existance}. This collection of solutions does not necessarily contains a unique solution. In this section we explore the properties of the extended difference operation and solutions collection $\mathcal{M}$.

%%%%%   Basic Properties   %%%%%
\subsection{Basic Properties}\label{subsec:basicProperties}
Corollary \ref{cor:properties1} shows that the extended difference satisfies certain desirable properties. For simplicity, when $A \ominus_e B$ is a singleton, we write $A \ominus_eB=X$ rather than $A \ominus_e B=\{X\}$.  
\begin{proposition}[Basic Properties]\label{cor:properties1}
Let \(A, B, C \in \KCn\). The extended set difference collection defined in (\ref{def:extended-diff}) satisfies the following properties:
\begin{enumerate}
    \item \(A \ominus_e A = \{0\}\).
    \item \((A \oplus B) \ominus_e B = A\).
    \item For $\lambda \in \mathbb{R}_{+}$, \(\lambda(A \ominus_e B) = \lambda A \ominus_e \lambda B\)
    \item For any $v \in \Rd$, $A\ominus_e (B\oplus\{v\}) = (A\ominus_eB) \oplus \{-v\} $
    \item For any $u \in \Rd$, $(A\oplus\{v\})\ominus_e B = (A \ominus_e B) \oplus \{v\}$.
    \item For $C\in \KCn$, $\left(A \ominus_e B\right)
    =(A \oplus C) \ominus_e(B \oplus C)$ 
\end{enumerate}
\end{proposition}
\begin{proof} \\
1. $d_H(A,A+\{0\})=0$, therefore $A \ominus_e A$ is unique and includes only $\{0\}$.\\
2. $d_H(A\oplus B,B\oplus A)=0$ and thus only $A$ can be in the difference set.\\
3. $X$ is an element of the extended distance $A \ominus_e B$ $\iff$ for any $X' \in \KCn$ and $\lambda \geq 0$, $\lambda d_H(A,B \oplus X) \leq \lambda d_H(A,B \oplus X')$ $\iff$ $d_H(\lambda A, \lambda (B \oplus X)) \leq d_H(\lambda A, \lambda (B \oplus X))$ $\iff$ $\lambda X$ is an element of $\lambda A \ominus_e \lambda B$.\\ 
4. Let $v \in \Rd$. $X$ is in $A\ominus_e B$ $\iff$ $d_H(A,B\oplus X)=d_H(A,(B\oplus\{v\})\oplus (X\oplus\{-v\})) \leq d_H(A,(B\oplus\{v\})\oplus (X' \oplus\{-v\}))$ for any $X' \in \KCn$ $\iff$ $X \oplus \{-v\}$ is in $A \ominus_e (B+\{v\})$.\\
5. The claim is proved in a similar way to (4).\\ 
6. $X$ is in $(A\oplus C) \ominus_e (B\oplus C)$ $\iff$ $d_H(A \oplus C,(B \oplus C) \oplus X) \leq d_H(A \oplus C,(B \oplus C) \oplus X')$, for $\forall X' \in \KCn$ $\iff$  $d_H(A ,B \oplus X) \leq d_H(A,B \oplus X')$, for $\forall X' \in \KCn$, $\iff$ $X$ is in $A \ominus_e B$. 
\end{proof}

Note: If $X_1, X_2 \in A \ominus_e B$, it is not necessarily true that $X_1 \cap X_2 \neq \emptyset$. It is also not necessarily true that $B \oplus X_1 \cap B \oplus X_2 \neq \emptyset$. An example in $\mathbb{R}^2$ is $A$ being the interval between $(-\frac{1}{2},0)$ and $(\frac{1}{2},0)$, and $B$ being the interval between $(0,-1)$ and $(0,1)$. In this case $d_H(A,B) = 1 \leq d_H(A,B+X), \forall X \in \KCn$. It is possible to show that $X_1= \{(-\frac{1}{2},0)\}$ and $X_2= \{(\frac{1}{2},0)\}$ both are in $A \ominus_e B$ and yet $X_1 \cap X_2 = \emptyset$ and $B \oplus X_1 \cap B \oplus X_2 = \emptyset$.

%%%%%    Convexity   %%%%%%%
\begin{proposition}
    The set of solutions, $\mathcal{M}$, is a convex set with respect to the Minkowski addition. Let $X_1, X_2 \in \mathcal{M}$, then for $\lambda \in [0,1]$, $\lambda X_1 + (1-\lambda) X_2 \in \mathcal{M}$.
\end{proposition}
\begin{proof}
    Let $d_H(A,B\oplus X_1)=d_H(A,B\oplus X_2)=\delta$. Then, $d_H(A,B \oplus \lambda X_1 + (1-\lambda) X_2) \leq \lambda d_H(A,B \oplus X_1) + (1 - \lambda)d_H(A,B \oplus X_2) = \delta$. Therefore, $\lambda X_1 + (1-\lambda) X_2 \in \mathcal{M}$.
\end{proof}

%%%%%    Transformations    %%%%%
\subsection{Orthonormal Transformations}
Consider an orthogonal transformations of sets in $\KCn$ as these transformations preserve length, convexity, and compactness. The following proposition shows that the order of orthogonal transformations and of extended set difference of sets in $\KCn$ are exchangeable.

\begin{proposition}\label{prop:orthogonalTransformations}
Let $T: \Rd \rightarrow \Rd$ be an orthogonal transformation, i.e, $T^{\top} T=I$, where $I$ is the identity matrix. For any,  $A, B \subset \KCn$, we have:

$$
T\left(A \ominus_e B\right)=T A \ominus_e T B,
$$
where $T(A \ominus_e B)$ means the collection of $TX$ for all $X \in A \ominus_e B$.
\end{proposition}
\begin{proof}
Let $X^*$ be in $A \ominus_e B$. $T(B \oplus X^*)=T B \oplus T X^*$, using the distributive property of $T$ over Minkowski summation. Since $T$ preserves the Hausdorff distance, we have:
$$
d_H(T A, T B \oplus T X^*)=d_H(A, B \oplus X^*).
$$
Therefore, since $X^*$ minimizes $d_H(A, B \oplus X)$, $T X^*$ minimizes $d_H(T A, T B \oplus T X)$. Hence, $T X^* \in T A \ominus_e T B$. The other direction of inclusion follows the same steps.
\end{proof}

%%%%   Convergence   %%%% 
\subsection{Convergence}
In this section, we show that the extended set difference and set limits are interchangeable. The convergence notion we use is convergence in the Hausdorff sense. For a sequence $\{A_n\}$ of elements of $\KCn$ and a set $A \in \KCn$, we say that $A_n \xrightarrow{d_H} A$ if $d_H(A_n,A) \rightarrow 0$ as $n \rightarrow \infty$ (see \cite{tuzhilin-2020}).

\begin{proposition}[Convergence]\label{Convergence}
Suppose $A_1,A_2,...$ and $B_1,B_2,... $ are two sequences of sets in $\KCn$, and $A$ and $B$ are two sets in $\KCn$ such that $A_n \xrightarrow{d_H} A$ and $B_n \xrightarrow{d_H} B$ in the Hausdorff sense. Then, there is a sequence $\{X_n\}$ such that $X_n \in A_n \ominus_e B_n$ and $X_n \xrightarrow{d_H} X^*$ and $X^*$ is in $A \ominus_e B$.
\end{proposition}

\begin{proof}
By the Blaschke Selection Theorem, the family $\KCn$ of compact convex sets in $\mathbb{R}^d$ is sequentially compact under the Hausdorff metric, so the sequence $\{X_n\}$ has a convergent subsequence $X_{n_k}\to X'$ for some $X'\in\KCn$. Since $A_{n_k} \xrightarrow{d_H} A$ and $B_{n_k}\xrightarrow{d_H} B$ and since Minkowski addition is 1-Lipschitz (thus continuous) with respect to $d_H$, we also have $B_{n_k}\oplus X_{n_k} \xrightarrow{d_H} B\oplus X'$ and $B_{n_k}\oplus X^*\xrightarrow{d_H} B\oplus X^*$, for any $X^*$ in $A\ominus_e B$. Consequently, 
\[
  d_H\bigl(A_{n_k},B_{n_k}\oplus X_{n_k}\bigr)
  \;\to\; 
  d_H\bigl(A,B\oplus X'\bigr)
  \quad\text{and}\quad
  d_H\bigl(A_{n_k},B_{n_k}\oplus X^*\bigr)
  \;\to\; 
  d_H\bigl(A,B\oplus X^*\bigr).
\]
By definition, each $X_{n_k}$ minimizes $d_H(A_{n_k}, B_{n_k}\oplus X)$ over $X\in\KCn$, so 
\[
  d_H\bigl(A_{n_k},B_{n_k}\oplus X_{n_k}\bigr)
  \;\le\; 
  d_H\bigl(A_{n_k},B_{n_k}\oplus X^*\bigr).
\]
Taking the limit as $k\to\infty$ shows 
\[
  d_H\bigl(A,B\oplus X'\bigr)
  \;\le\;
  d_H\bigl(A,B\oplus X^*\bigr).
\]
On the other hand, $X^*$ itself is defined to minimize $d_H(A,B\oplus X)$, so $d_H(A,B\oplus X^*) \le d_H(A,B\oplus X')$ holds for all $X'\in\KCn$. Hence
\[
  d_H\bigl(A,B\oplus X^*\bigr)
  \;\le\; 
  d_H\bigl(A,B\oplus X'\bigr)
  \;\le\;   d_H\bigl(A,B\oplus X^*\bigr),
\]
and therefore $d_H(A,B\oplus X') = d_H(A,B\oplus X^*)$. By the definition of $A \ominus_e B$, it follows that $X'$ is in $A \ominus_e B$. Since $X'$ is the limit of an arbitrary convergent subsequence of $\{X_n\}$, every convergent subsequence must have the same limit $X'$, which implies that the entire sequence $X_n$ converges to $X'$.
 \end{proof}

%%%%%    scope   %%%%%
\subsection{Scope}\label{subsec:scope}

\begin{theorem}[Scope]\label{thm:Scope}
Let $A, B \in \KCn$. Let $X \in A \ominus_e B$ and let $\delta=d(A,B \oplus X)$ be the minimum distance. Then, for every two elements $X_1, X_2 \in A \ominus_e B$, $d(X_1,X_2) \leq 2\delta$. Moreover, for any $K \in \KCn$, $d(X_1,X_2) \leq 2d(A,B \oplus K)$ including for $K=\{0\}$. Moreover, if $\exists X \in A \ominus_e B$ such that $A = B \oplus X$, then $A \ominus_e B$ is a singleton.
\end{theorem}

\begin{proof}
Assume \(X_1, X_2 \in \KCn\) are such that $d_H(A,B \oplus X_1)=d_H(A,B \oplus X_2)=\delta$. By triangle inequality,  
\begin{align*}
    d(X_1,X_2) &\leq d_H(B \oplus X_1, B \oplus X_2) \\
               &\leq d_H(B \oplus X_1, A) + d_H(A, B \oplus X_2) \\
               &=2 \delta \\
               &\leq 2d_H(A,B \oplus K),
\end{align*}
where the last inequality comes from the fact that $K \notin A \ominus_e B$. Finally, if $\exists X$ such that $A=B \oplus X$, then $\delta=0$. For $X' \in A \ominus_e B$, it must be $d_H(A,B \oplus X')=\delta=0$ and thus, $A=B \oplus X'$. Then, $B \oplus X=B \oplus X'$. Since all these sets are compact, $X=X'$.
\end{proof}

For a set $K\in \KCn$ let $\rin(K)$ be the radius of the largest inscribed ball in $K$ and let $\Rout(K)$ be the radius of the smallest ball that contains $K$.

\begin{theorem}\label{thm:radiusBounds}
    Let $X \in A \ominus_e B$, then
    \begin{enumerate}
        \item $max\{0,\rin(A)-\Rout(B)-\delta\} \leq \rin(X) \leq max\{0,\rin(A)-\rin(B)+\delta\}$.
        \item max\{$\Rout(A)-\Rout(B)-\delta\} \leq \Rout(X) \leq max\{0,\Rout(A)-\rin(B)+\delta\}.$
    \end{enumerate}
\end{theorem}

\begin{proof}
    Let $\mathbb{B}$ denote the unit ball centered at the origin and $K, L \in \KCn$. We start by stating the following properties for any $K,L \in \KCn$.
    (a) if $K\subset L$, $\rin(K)\leq \rin(L)$ and $\Rout(K) \leq \Rout(L)$, (b) For any $\delta \geq 0$, $\rin(K \oplus \delta \mathbb{B})=\rin(K)+\delta$ and $\Rout(K \oplus \mathbb{B}_t)=\Rout(K)+\delta$, (c) $d_H(K,L)=\delta$ implies $K \subset L \oplus \delta \mathbb{B}$ and $L \subset L \oplus \delta \mathbb{B}$. (d) $K \subset \Rout(K)\mathbb{B}$. These properties imply that (e) $d_H(K,L)=\delta$ implies $\rin(K)\leq \rin(K)+\Rout(L) +\delta$
    \begin{enumerate}     
    \item By (e) $\rin(A) \leq \rin(X) + \Rout(B) + \delta$. This inequality and the non-negativity of $\rin(X)$ imply the lower bound. Since $B \oplus X \subset A + \delta \mathbb{B}$, by (a) and (b) $\rin(B\oplus X) \leq \rin(A) +\delta$. Moreover, $\rin(B)+\rin(X) \leq \rin(B\oplus X)$. Combining these two inequalities implies the upper bound.
        \item The same as (1) above by replacing $\rin()$ with $\Rout()$.
    \end{enumerate}
\end{proof} %This proof is now done

\textbf{Note:} If $\exists X$ such that $A=B\oplus X$, then $\delta=0$ and the bounds in Therom \ref{thm:radiusBounds} tighten.
%%%%   Uniqueness Refinement   %%%%

\begin{example}
    Consider $A=[-1,1] \times \{0\}$ and $B=\{0\}\times [-1,1]$. The minimizer set is the family $X=[-t,t]\times \{0\}$ for $t\in[0,1]$. $\rin(A)=\rin(B)=0$ and $\Rout(A)=\Rout(B)=1$. Let $X_1=\{-1\}\times[-1,1]$ and $X_2=\{1\}\times[-1,1]$. Then, $\delta = d_H(A,B)=1$ and $d_H(X_1,X_2)=2$ satisfies the inequality in Theorem \ref{thm:Scope} with equality. Moreover, the first inequality in Theorem \ref{thm:radiusBounds} is satisfied with equality for the lower bound and strong inequality for the upper bound for both $\rin(X_1)$ and $\rin(X_2)$. For the second set of inequalities the lower bound is satisfied weakly and the upper bound with equality for both $\Rout(X_1)$ and $\Rout(X_2)$.  
\end{example}

%%% Symmetry %%%
\subsection{Symmetry}
Let $\mathcal{O}(d)$ be the group of orthonormal transformations in $\Rd$. Let $G \subset \mathcal{O}(d)$ be a finite subgroup of orthogonal transformations.
\begin{definition}
    We say that $A \in \KCn$ is $G$-symmetric about a point $c \in \Rd$, if $A$ is invariant for $g(x) = c+\!A_g(x-c)$ in $G$.
\end{definition}

Next we show that if $A$ is $G$-symmetric about $c_A$ and $B$ is $G$-symmetric about $c_B$ (same linear action $g(x) = c+\!A_g(x-c)$ with $A_g\in \mathcal{O}(d)$), then $A\ominus_e B$ always contains minimizers that are $G$-symmetric about $c_A-c_B$.

Let

$$
\delta=\inf_Y d_H(A,\,B\oplus Y),\qquad
\mathcal M\ :=\ \bigl\{X\in\KCn:\ d_H\!\bigl(A,\ B\oplus X\bigr)=\delta\bigr\}.
$$

\begin{theorem}
    Assume $g(A)=A$ via $g(x)=c_A+A_g(x-c_A)$ and $g(B)=B$ via $g(x)=c_B+A_g(x-c_B)$ for all $g\in G$.

Then:

    1. $G$-invariance of the minimizer set (about $c_A-c_B$).
   For every $X\in\mathcal M$ and $g\in G$, the set

$$
g\!\cdot\!X\;:=\;(c_A-c_B)\;+\;A_g\big(X-(c_A-c_B)\big)
$$

    is also in $\mathcal M$.

    2. Existence of a $G$-symmetric minimizer.
   There exists $\overline X\in\mathcal M$ such that $g\!\cdot\!\overline X=\overline X$ for all $g\in G$; i.e., $\overline X$ is $G$-symmetric about $c_A-c_B$.

    3. Uniqueness implies symmetry.
   If $\mathcal M$ is a singleton, its (unique) member is $G$-symmetric about $c_A-c_B$.

\end{theorem}
\begin{proof}
Let $G\subset\mathcal O(d)$ be finite and suppose $A,B\in\KCn$ are $G$-symmetric about $c_A,c_B$, respectively. Therefore,
$$
g(x)=c_A+A_g(x-c_A)\ \text{ maps }A\text{ to itself, and }\ 
g(x)=c_B+A_g(x-c_B)\ \text{ maps }B\text{ to itself for all }g\in G,
$$
with $A_g\in\mathcal O(d)$. Define the affine $G$-action on sets about $c_A-c_B$ by
\[
g\!\cdot\!X\ :=\ (c_A-c_B)\,+\,A_g\bigl(X-(c_A-c_B)\bigr).
\]

\smallskip
Let,
\[
A':=A-c_A,\qquad B':=B-c_B.
\]
For any $X\in\KCn$ set
\[
Y:=X-(c_A-c_B).
\]
Since Hausdorff distance is invariant under a common translation and Minkowski sum commutes with translation, we have
\begin{align*}
d_H\!\bigl(A,\ B\oplus X\bigr)
&= d_H\!\bigl(A-c_A,\ (B\oplus X)-c_A\bigr) \\
&= d_H\!\bigl(A',\ (B-c_B)\oplus (X-(c_A-c_B))\bigr)\\
&= d_H\!\bigl(A',\ B'\oplus Y\bigr).
\end{align*}
Thus $X\in\mathcal M$ if and only if $Y$ minimizes $f(Y):=d_H(A',B'\oplus Y)$. Moreover, by the assumption on $A,B$, the normalized sets $A',B'$ are invariant under the \emph{linear} action $x\mapsto A_gx$ for every $g\in G$.

%\smallskip
%\noindent\emph{Support-function representation and orthogonal invariances.} Recall that for compact convex $K,L$,
%\[
%d_H(K,L)\ =\ \sup_{\|u\|=1}\bigl|h_K(u)-h_L(u)\bigr|,\qquad
%h_{K\oplus L}=h_K+h_L,\qquad h_{A_gK}(u)=h_K(A_g^\top u).
%\]
%Hence orthogonal maps preserve $d_H$:
%\[
%d_H(A_gK,A_gL)=d_H(K,L),
%\]
%and commute with Minkowski sums:
%\[
%A_g(K\oplus L)=A_gK\ \oplus\ A_gL.
%\]

\smallskip
\noindent\textbf{(1) $G$-invariance of the minimizer set.}
Let $Y$ be a minimizer of $f$, and fix $g\in G$. Then
\begin{align*}
f(A_gY)
&= d_H\!\bigl(A',\, B'\oplus A_gY\bigr)
= d_H\!\bigl(A',\, A_g(B'\oplus Y)\bigr)\\
&= d_H\!\bigl(A_gA',\, A_g(B'\oplus Y)\bigr)
= d_H\!\bigl(A',\, B'\oplus Y\bigr)
= f(Y),
\end{align*}
using that $A_gA'=A'$ and $A_gB'=B'$.
Thus $A_gY$ is also a minimizer. Translating back, if $X\in\mathcal M$ then
\[
g\!\cdot\!X\ =\ (c_A-c_B)+A_g\bigl(X-(c_A-c_B)\bigr)\ \in\ \mathcal M.
\]

\smallskip
\noindent\textbf{(2) Existence of a $G$-symmetric minimizer.}
Pick any minimizer $Y$ of $f$. Define the Minkowski average of its $G$-orbit by
\[
\overline{Y}\ :=\ \frac{1}{|G|}\ \bigoplus_{g\in G}\ A_gY,
\]
equivalently for support functions,
\[
h_{\overline{Y}}(u)\ :=\ \frac{1}{|G|}\sum_{g\in G} h_{A_gY}(u),\qquad \|u\|=1.
\]
Write $\phi(u):=h_{A'}(u)-h_{B'}(u)$. Then, for the objective
\[
f(Z)=\|\phi-h_Z\|_\infty,
\]
the convexity of the sup-norm yields
\begin{align*}
f(\overline{Y})
&= \left\|\phi-\frac{1}{|G|}\sum_{g} h_{A_gY}\right\|_\infty
\ \le\ \frac{1}{|G|}\sum_{g}\|\phi-h_{A_gY}\|_\infty \\
 &= \frac{1}{|G|}\sum_{g} f(A_gY)
\ =\ f(Y).
\end{align*}
Since $Y$ is optimal, equality holds; hence $\overline{Y}$ is a minimizer.

We now show $\overline{Y}$ is $G$-invariant. For any $h\in G$ and any unit $u$,
\begin{align*}
h_{A_h\overline{Y}}(u)
&= h_{\overline{Y}}(A_h^\top u)
= \frac{1}{|G|}\sum_{g\in G} h_{A_gY}(A_h^\top u)
= \frac{1}{|G|}\sum_{g\in G} h_{A_hA_gY}(u)\\
&= \frac{1}{|G|}\sum_{g'\in G} h_{A_{g'}Y}(u)
= h_{\overline{Y}}(u),
\end{align*}
so $A_h\overline{Y}=\overline{Y}$ for all $h\in G$. Translating back, $\overline{X}:=\overline{Y}+(c_A-c_B)\in\mathcal M$ and
\[
g\!\cdot\!\overline{X}=(c_A-c_B)+A_g\bigl(\overline{X}-(c_A-c_B)\bigr)=\overline{X},
\]
i.e., $\overline{X}$ is $G$-symmetric about $c_A-c_B$.

\smallskip
\noindent\textbf{(3) Uniqueness implies symmetry.}
If $\mathcal M=\{X^\star\}$ is a singleton, part (1) gives $g\!\cdot\!X^\star\in\mathcal M$ for every $g\in G$, hence $g\!\cdot\!X^\star=X^\star$. Thus $X^\star$ is $G$-symmetric about $c_A-c_B$.
%
%\smallskip
%All three claims follow.
\end{proof}

Therefore, when $A$ and $B$ are symmetric about different centers, their extended difference always contains minimizers symmetric about the difference of the centers $c_A-c_B$.

%%%%%%%%%%%%%%%%%%%%%%%%%%%%%%%%%%%%%%%%%%%%%%%%%
%%%%%    Uniqueness through Refinement     %%%%%%
%%%%%%%%%%%%%%%%%%%%%%%%%%%%%%%%%%%%%%%%%%%%%%%%%
\section{Uniqueness Refinement}\label{sec: refinement}
To achieve a unique solution in the optimization problem in Definition \ref{def:extended-diff}, we add a perturbation term, ensuring the objective functional is strictly convex. Before doing so, Lemma  \ref{lemm:convexoriginal}, shows that the original objective functional proposed in Definition \ref{def:extended-diff} is convex.  

\begin{lemma}\label{lemm:convexoriginal}
Let \(A,B\in \KCn\) be fixed and for each \(X\in \KCn\), define
\[
f(X) = d_H\bigl(A,\,B \oplus X\bigr).
\]

The functional $f(X)$ is convex. That is, for any two sets \(X_1,X_2\in \KCn\) and any \(\lambda\in[0,1]\), if we define
\[
X_\lambda = \lambda X_1\oplus (1-\lambda)X_2,
\]
we have
\[
f(X_\lambda) \le \lambda f(X_1) + (1-\lambda)f(X_2).
\]
\end{lemma}

\begin{proof}
Note that since Minkowski addition is linear with respect to the support function, we have for any \(X\in \KCn\)
\[
h_{B\oplus X}(u)=h_B(u)+h_X(u) \quad \text{for all } u\in\mathbb{R}^n.
\]
Hence,
\[
f(X)= d_H\bigl(A,\,B\oplus X\bigr)
=\sup_{\|u\|=1}\Bigl| h_A(u)-\bigl[h_B(u)+h_X(u)\bigr]\Bigr|.
\]
Defining
\[
\phi(u)=h_A(u)-h_B(u),
\]
we rewrite the above as
\[
f(X)=\sup_{\|u\|=1}\Bigl|\phi(u)-h_X(u)\Bigr|.
\]

For any two sets \(X_1,X_2\in \KCn\) and any \(\lambda\in[0,1]\), let
\[
X_\lambda=\lambda X_1\oplus (1-\lambda)X_2.
\]
The linearity of the support function under Minkowski sums gives
\[
h_{X_\lambda}(u)=\lambda h_{X_1}(u)+(1-\lambda)h_{X_2}(u).
\]
Thus, for all \(u\) with \(\|u\|=1\),
\[
\Bigl|\phi(u)-h_{X_\lambda}(u)\Bigr|
=\Bigl|\phi(u)-\lambda h_{X_1}(u)-(1-\lambda)h_{X_2}(u)\Bigr|.
\]
The convexity of the absolute value then implies
\[
\Bigl|\phi(u)-h_{X_\lambda}(u)\Bigr|
\le \lambda\Bigl|\phi(u)-h_{X_1}(u)\Bigr|+(1-\lambda)\Bigl|\phi(u)-h_{X_2}(u)\Bigr|.
\]
Taking the supremum over \(u \in \Sd\) yields
\[
\begin{aligned}
f(X_\lambda)&=\sup_{u \in \Sd}\Bigl|\phi(u)-h_{X_\lambda}(u)\Bigr|
\\[1mm]
&\le \lambda\,\sup_{u \in \Sd}\Bigl|\phi(u)-h_{X_1}(u)\Bigr|+(1-\lambda)\,\sup_{u \in \Sd}\Bigl|\phi(u)-h_{X_2}(u)\Bigr|
\\[1mm]
&=\lambda f(X_1)+(1-\lambda)f(X_2).
\end{aligned}
\]
This completes the proof.
\end{proof}

\begin{proposition}  
[Uniqueness via Strictly Convex Perturbation]\label{prop:general-uniqueness}
Let \( A, B \in \KCn \), and assume \( A \ominus_e B = \arg \inf_{X \in \KCn} d_H(A, B \oplus X) \) is non-empty, as established by Theorem \ref{thm:Existance}. Consider the perturbed problem \( \arg \min_{X \in \KCn} \left[ d_H(A, B \oplus X) + \gamma R(X) \right] \), where \( R : \KCn \to [0, \infty) \) is a strictly convex, lower semicontinuous functional, and \( \gamma > 0 \) is a small positive constant. The perturbed objective functional ensures a unique solution for any \( A, B \in \KCn \).
\end{proposition}
\begin{proof}
Theorem \ref{thm:Existance} guarantees that \( A \ominus_e B = \arg \inf_{X \in \KCn} d_H(A, B \oplus X) \) is non-empty. For a strictly convex function $R:\KCn \rightarrow [0,\infty)$, define \( f_\gamma(X) = d_H(A, B \oplus X) + \gamma R(X) \). The addition of the strictly convex term \( \gamma R(X) \) to the convex \( d_H(A, B \oplus X) \) makes \( f_\gamma \) strictly convex on \(\KCn\), ensuring the optimization problem has a unique solution \( X^* \) for any \( A, B \in \KCn \).
\end{proof}

\begin{example}[Example of a Strictly Convex Perturbation]
A suitable choice for the strictly convex functional is \( R(X) = \int_{\Sd} h_{X}(u)^2 \, du \), where the integral is taken over the unit sphere \( \Sd \) in \( \Rd \).  To verify strict convexity, consider \( X_1, X_2 \in \KCn \) with \( X_1 \neq X_2 \), and let \( X_{\lambda} = \lambda X_1 + (1 - \lambda) X_2 \) for \( 0 < \lambda < 1 \). The support function of the Minkowski sum satisfies \( h_{X_{\lambda}}(u) = \lambda h_{X_1}(u) + (1 - \lambda) h_{X_2}(u) \), as \( \oplus \) is linear. Define \( R(X) = \int_{\Sd} h_{X}(u)^2 \, du \). Then,
\[
R(X_{\lambda}) = \int_{\Sd} \left[ \lambda h_{X_1}(u) + (1 - \lambda) h_{X_2}(u) \right]^2 \, du.
\]
Expanding the integrand, 
\[
R(X_{\lambda}) = \lambda^2 \int_{\Sd} h_{X_1}(u)^2 \, du + 2\lambda(1 - \lambda) \int_{\Sd} h_{X_1}(u) h_{X_2}(u) \, du + (1 - \lambda)^2 \int_{\Sd} h_{X_2}(u)^2 \, du.
\]
By comparison, the convex combination is:
\[
\lambda R(X_1) + (1 - \lambda) R(X_2) = \lambda \int_{\Sd} h_{X_1}(u)^2 \, du + (1 - \lambda) \int_{\Sd} h_{X_2}(u)^2 \, du.
\]
Taking the difference,
\[
\lambda R(X_1) + (1 - \lambda) R(X_2) - R(X_{\lambda}) = \lambda(1 - \lambda) \int_{\Sd} \left[ h_{X_1}(u)^2 + h_{X_2}(u)^2 - 2 h_{X_1}(u) h_{X_2}(u) \right] \, du.
\] 
Since \( h_{X_1}(u)^2 + h_{X_2}(u)^2 - 2 h_{X_1}(u) h_{X_2}(u) = \left[ h_{X_1}(u) - h_{X_2}(u) \right]^2 \), we have:
\[
\lambda R(X_1) + (1 - \lambda) R(X_2) - R(X) = \lambda(1 - \lambda) \int_{\Sd} \left[ h_{X_1}(u) - h_{X_2}(u) \right]^2 \, du.
\]
For \( X_1 \neq X_2 \), \( h_{X_1} \neq h_{X_2} \) (as support functions uniquely determine compact convex sets), and since \( \left[ h_{X_1}(u) - h_{X_2}(u) \right]^2 \geq 0 \) with strict inequality on a set of positive measure on \( \Sd \), the integral is positive. Thus, for \( 0 < \lambda < 1 \),
\[
R(X) < \lambda R(X_1) + (1 - \lambda) R(X_2),
\]
confirming that \( R(X) \) is strictly convex on \( \KCn \).
\end{example}
The following proposition ensures that the modified problem defined in Proposition \ref{prop:general-uniqueness} yields a solution which approaches an optimal set in the original problem defined in Definition \ref{def:extended-diff}.
\begin{proposition}[Convergence of Perturbed Solution]\label{prop:convpert}
By Theorem \ref{thm:Existance}, the original problem in (\ref{def:extended-diff}) has a non-empty set of optimal solutions \( \mathcal{M} \). Let 
\begin{equation}
 X^*_{\gamma} = \arg \min_{X \in \KCn} [d_H(A, B \oplus X) + \gamma R(X)]  
\end{equation}
 where \( R: \KCn \to [0, \infty) \) is strictly convex and lower semicontinuous and \( \gamma > 0 \). As \( \gamma \to 0^+ \), \( X^*_{\gamma} \) converges in the Hausdorff metric to one of the solutions in \( \mathcal{M} \).
\end{proposition}

\begin{proof}
Let \( X^*_{\gamma} = \arg \min_{X \in \KCn} f_\gamma(X) \), where \( f_\gamma(X) = d_H(A, B \oplus X) + \gamma R(X) \), over the compact set \( \KCn | \mathbb{B}_r \) with \( A \subset B \oplus \mathbb{B}_r \). \( X^*_{\gamma} \) is unique due to the strict convexity of \( f_\gamma \). Suppose \( X^*_{\gamma} \) does not converge to any point in \( \mathcal{M} \). For \( \gamma_n = 1/n \), compactness of \( \KCn | \mathbb{B}_r \) implies a subsequence \( X^*_{\gamma_n} \to X_0 \in \KCn | \mathbb{B}_r \). For \( X \in \mathcal{M} \), \( d_H(A, B \oplus X) = \delta \), and:
\[
d_H(A, B \oplus X^*_{\gamma_n}) + \gamma_n R(X^*_{\gamma_n}) \leq \delta + \gamma_n R(X),
\]
so \( d_H(A, B \oplus X_{\gamma_n}) \leq \delta + \gamma_n [R(X) - R(X_{\gamma_n})] \). With \( R \) bounded by \( M \), \( \gamma_n [R(X) - R(X_{\gamma_n})] \leq 2 \gamma_n M \to 0 \), and since \( d_H \) is continuous, \( d_H(A, B \oplus X_{\gamma_n}) \to d_H(A, B \oplus X_0) \), thus \( d_H(A, B \oplus X_0) \leq \delta \), hence \( d_H(A, B \oplus X_0) = \delta \) and \( X_0 \in \mathcal{M} \). 

We now investigate how different choices of the strictly convex perturbation functional affect the limiting solution, particularly when the collection contains multiple elements.

Now, \( \delta \leq f_{\gamma_n}(X_{\gamma_n}) \leq \delta + \gamma_n R(X) \), so \( f_{\gamma_n}(X_{\gamma_n}) \to \delta \), implying \( d_H(A, B \oplus X_{\gamma_n}) \to \delta \) and \( \gamma_n R(X_{\gamma_n}) \to 0 \). If \( X_\gamma \) does not approach \( \mathcal{M} \), there exists \( \eta > 0 \) with \( d_H(X_{\gamma_n}, \mathcal{M}) \geq \eta \) for some sequence, but all limit points are in \( \mathcal{M} \), a contradiction. Thus, \( X_\gamma \to X_0 \in \mathcal{M} \), where \( X_0 \) is unique in \( \mathcal{M} \) due to \( R \)’s strict convexity. Hence, \( X^* \to X_0 \in \mathcal{M} \) as \( \gamma \to 0^+ \).
\end{proof}
 We now demonstrate that distinct choices of the strictly convex perturbation functional generally yield different limiting solutions when multiple minimizers exist.

\begin{remark}\label{prop:distinct-limits}
Let \( A, B \in \KCn \), and let \( \mathcal{M} = A \ominus_e B = \arg \inf_{X \in \KCn} d_H(A, B \oplus X) \), with \( \delta = \inf_{X \in \KCn} d_H(A, B \oplus X) \). Consider two distinct strictly convex, lower semicontinuous functionals \( R, R' : \KCn \to [0, \infty) \), and define the perturbed objective functions for \( \gamma > 0 \):
\[
f_\gamma(X) = d_H(A, B \oplus X) + \gamma R(X), \quad f'_\gamma(X) = d_H(A, B \oplus X) + \gamma R'(X).
\]
Let \( X^*_\gamma = \arg \min_{X \in \KCn} f_\gamma(X) \) and \( X^{'*}_\gamma = \arg \min_{X \in \KCn} f'_\gamma(X) \). Assuming \( X^*_\gamma \to X^* \in \mathcal{M} \) and \( X^{'*}_\gamma \to X^{'*} \in \mathcal{M} \) as \( \gamma \to 0^+ \), where \( X^* = \arg \min_{X \in \mathcal{M}} R(X) \) and \( X^{'*} = \arg \min_{X \in \mathcal{M}} R'(X) \), if \( \mathcal{M} \) contains more than one element and \( R \neq R' \), then generally \( X^* \neq X^{'*} \).

To show this, suppose \( \mathcal{M} \) contains at least two distinct elements, say \( X_1, X_2 \in \mathcal{M} \), with \( X_1 \neq X_2 \), such that \( d_H(A, B \oplus X_1) = d_H(A, B \oplus X_2) = \delta \). Since \( R \) and \( R' \) are strictly convex and lower semicontinuous, their minimizers over the compact, convex set \( \mathcal{M} \):
\[
X^* = \arg \min_{X \in \mathcal{M}} R(X), \quad X^{'*} = \arg \min_{X \in \mathcal{M}} R'(X),
\]
are unique due to strict convexity. We prove that \( X^* \neq X^{'*} \) in the general case by showing that distinct functionals typically select different elements from \( \mathcal{M} \).

Define:
\[
R(X) = \int_{\mathbb{S}^{d-1}} [h_X(u) - h_{X_1}(u)]^2 \, du, \quad R'(X) = \int_{\mathbb{S}^{d-1}} [h_X(u) - h_{X_2}(u)]^2 \, du.
\]
As established in Subsection 2.3, $R(X)$ and $R'(X)$ are strictly convex and lower semicontinuous. Evaluate:
- For \( X = X_1 \), \( h_{X_1}(u) - h_{X_1}(u) = 0 \), so \( R(X_1) = 0 \). Since \( X_1 \neq X_2 \), the support functions differ on a set of positive measure (as \( h_X \) uniquely determines \( X \in \KCn \)), so \( R(X_2) = \int_{\mathbb{S}^{d-1}} [h_{X_2}(u) - h_{X_1}(u)]^2 \, du > 0 \). Thus, \( X^* = X_1 \). For \( X = X_2 \), \( h_{X_2}(u) - h_{X_2}(u) = 0 \) and \( h_{X_2}(u) - h_{X_1}(u) \) are distinct, so \( R'(X_2) = 0 \) and \( R'(X_1) > 0 \). Thus, \( X^{'*} = X_2 \). Since \( X_1 \neq X_2 \), it follows that \( X^* \neq X^{'*} \).
To generalize, since \( R \neq R' \), their strict convexity ensures unique minimizers over \( \mathcal{M} \). When \( \mathcal{M} \) contains multiple elements, the minimizers \( X^* \) and \( X^{'*} \) differ unless \( R \) and \( R' \) coincidentally share the same minimizer, which requires \( R(X) = R'(X) + c \) for some constant \( c \) over \( \mathcal{M} \). Such alignment is non-generic, as \( R \) and \( R' \) are arbitrary distinct functionals. Thus, \( X^* \neq X^{'*} \) generally holds.
If \( \mathcal{M} = \{ X_0 \} \), then \( X^* = X^{'*} = X_0 \), as both functionals minimize at the unique element.
\end{remark}

\begin{example}\label{ex:distinct-limits}
Let,
\[
A = [-1, 1] \times \{0\}, \quad B = \{0\} \times [-1, 1],
\]
where \( A \) is the line segment from \( (-1, 0) \) to \( (1, 0) \), and \( B \) is the line segment from \( (0, -1) \) to \( (0, 1) \). Denote the extended set difference as \( \mathcal{M} = A \ominus_e B = \arg \inf_{X \in \KCn} d_H(A, B \oplus X) \).

For \( X = [-t, t] \times \{0\} \), \( t \geq 0 \), the Minkowski sum is $B \oplus X = [-t, t] \times [-1, 1]$ ,
a rectangle centered at the origin with width \( 2t \) and height 2. It is easy to verify that for $t \in [0,1]$, 
$d_H(A, B \oplus X) = 1$ and for $t>1$, $d_H(A, B \oplus X) > 1$.
Therefore, the infimum is achieved at any \( t \in [0, 1] \). Thus,
\[
\mathcal{M} = \{ [-t, t] \times \{0\} \mid t \in [0, 1] \},
\]
a non-singleton set including \( \{(0,0)\} \) (at \( t = 0 \)) and \( [-1, 1] \times \{0\} \) (at \( t = 1 \)).

Consider the penalized objective function:
\[
f_\gamma(X) = d_H(A, B \oplus X) + \gamma R(X),
\]
and test two strictly convex, lower semicontinuous penalty functionals:
\[
R(X) = \int_{\mathbb{S}^1} h_X(u)^2 \, du, \quad R'(X) = \int_{\mathbb{S}^1} [h_X(u) - h_{X_2}(u)]^2 \, du,
\]
where \( X_2 = [-1, 1] \times \{0\} \). For \( X = [-t, t] \times \{0\} \), \( h_X(u) = t |u_1| \), where \( u = (u_1, u_2) \). Compute:
\[
R(X) = \int_{\mathbb{S}^1} (t |u_1|)^2 \, du = t^2 \int_0^{2\pi} \cos^2 \theta \frac{d\theta}{2\pi} = t^2 \cdot \frac{1}{2} = \frac{t^2}{2}.
\]
For \( X_2 \), \( h_{X_2}(u) = |u_1| \), so:
\[
R'(X) = \int_{\mathbb{S}^1} (t |u_1| - |u_1|)^2 \, du = (t - 1)^2 \int_{\mathbb{S}^1} |u_1|^2 \, du = \frac{(t - 1)^2}{2}.
\]
Since \( d_H(A, B \oplus X) = 1 \) for \( X \in \mathcal{M} \), the limiting solution as \( \gamma \to 0^+ \) minimizes \( R(X) \) or \( R'(X) \) over \( \mathcal{M} \):
\begin{itemize}
\item For \( R(X) = \frac{t^2}{2} \), the minimum occurs at \( t = 0 \), so \( X^* = \{(0,0)\} \).
\item For \( R'(X) = \frac{(t - 1)^2}{2} \), the minimum occurs at \( t = 1 \), so \( X^{'*} = [-1, 1] \times \{0\} \).
\end{itemize}
Therefore, using \( R \) versus \( R' \) in the penalized objective function yields distinct limits \( X^* \neq X^{'*} \), confirming that different strictly convex penalty functionals generally select different elements from a non-singleton \( \mathcal{M} \).
\end{example}

%%%%%%%%%%%%%%%%%%%%%%%%%%%%%%%%%%%%%%%%%%
%%%  Computing the extended difference %%%
%%%%%%%%%%%%%%%%%%%%%%%%%%%%%%%%%%%%%%%%%%

\section{Computation of the Extended Set Difference}\label{sec:computation}

In this section, we propose an algorithm for finding an approximately optimal convex and compact set $X\in \KCn$ which is a minimizer of the Hausdorff distance between $A$ and $B\oplus X$ where $A,B, \text{and }X \in \KCn$. In Section \ref{sec:LPformulation} we show how to formulate the task of finding the extended difference as a Linear Programming (LP) problem. Although the method is applicable for any finite $d \geq 2$, in Section \ref{sec:computational} we demonstrate our algorithm for $d=2$. The examples in Section \ref{sec:Examples} demonstrate the algorithm.

We start by introducing some basic tools and results from convex analysis. Let \( C \subseteq \Rd \) be a nonempty, compact set. Then its support function \( h_C: \Rd \to \mathbb{R} \) is sublinear if and only if \( C \) is \textbf{convex} (\cite{hiriart2004fundamentals}, Theorem 3.1.1).\footnote{If $\mathcal{X}$ is a vector space, a function \(h : \mathcal{X} \to \mathbb{R}\) is said to be \emph{sublinear} if for every $x,y \in \mathcal{X}$, $h(x+y) \leq h(x) + h(y)$; and for every \(x \in \mathcal{X}\) and every scalar \(\lambda \ge 0\), $h(\lambda x) = \lambda\, h(x)$ (\cite{rockafellar-2015}).} Moreover, the space of sublinear functions on $\Rd$, denoted by \( \mathcal{S} \), forms a convex cone.

\par
Since every convex set $X$ has a sublinear support function $h_X$, and every sublinear function is the support function of a unique closed convex set (see Chapter 13, \cite{ rockafellar-2015}), we can convert our initial optimization problem (\ref{def:extended-diff}) to an optimization over the space of support functions. The following equivalence holds:

\begin{equation}\label{equ:LP}
  \min_{X \in \KCn} d_H(A, B + X) \equiv \min_{h \in \mathcal{S}} \|f - h\|_{\infty}
\end{equation}
where, $\mathcal{S}$ is the convex cone of sublinear functions on $\Sd$ and $f=h_A-h_B$. The cone property (i.e., closure under nonnegative scalar multiplication and addition) ensures that the set over which we minimize is convex which is essential for formulating (\ref{equ:LP}) as an LP problem.

\subsection{LP Formulation}\label{sec:LPformulation}
In this subsection, we describe a linear programming (LP) approach to solve the optimization problem in (\ref{equ:LP}). Our method relies on the following topological properties of $\Sd$: 

\begin{enumerate}
  \item \textbf{Compactness:} \\
  \( \Sd \) is a compact subset of \( \Rd \). As a result:
  \begin{itemize}
    \item Every continuous function \( f: \Sd \to \mathbb{R} \) (such as $ h_A$, $h_B$ or their difference) attains its maximum and minimum on $\Sd$.
    \item Every continuous function on \( \Sd \) is uniformly continuous. This implies that for any \( \epsilon > 0 \), there exists \( \delta > 0 \) such that for all \( u, v \in \Sd \),
    \[
    \|u - v\| < \delta \implies |f(u) - f(v)| < \epsilon.
    \]
  \end{itemize}

  \item \textbf{Separability (Polish Space Property):} \\
  \( \Sd \) is a Polish space, meaning it is a complete, separable metric space. As a result:
  \begin{itemize}
    \item There exists a countable dense subset \( \{u_1, u_2, \dots\} \subset \Sd \).
    \item For any \( \delta > 0 \), there exists a finite (or countable) set of directions \( \mathcal{U} = \{u_1, u_2, \dots, u_m\} \subset \Sd \) that is \( \delta \)-dense in \( \Sd \). That is, for every \( u \in \Sd \), there exists \( u_i \in \mathcal{U} \) such that
    \[
    \|u - u_i\| < \delta.
    \]
  \end{itemize}
\end{enumerate}
Building on these properties and since support functions are uniformly continuous on the compact unit sphere, we take the following approach. Instead of optimizing over all directions $u \in$ $\mathbb{S}^{n-1}$, we select a finite set of \(m\) points in $\Sd$ (for instance, in $\mathbb{R}^2$, by setting \(u_i=(\cos\theta_i,\sin\theta_i)\) for \(\theta_i=2\pi i/m,\ i=0,\dots,m-1\)), 
$$
\mathcal{U}=\left\{u_1, u_2, \ldots, u_m\right\} \subset \Sd.
$$
Let $f(u)=h_A(u)-h_B(u)$ and let 
\[
f_i=f(u_i),\quad i=1,\dots,m.
\]

Given the set of directions $\mathcal{U}$, we approximate the support function \(h_X\) only on $\mathcal{U}$. To obtain a best uniform (Chebyshev) approximation of $h_X$ we seek $m$ values, $x_0,...,x_{m-1}$, approximating $h_X(u_i)$, respectively, such that (1) the sub linearity of $h_X$ is maintained and (2) the distance to $f(u)=h_A(u)-h_B(u)$ is minimized. We consider the following linear program (LP):

\begin{equation}\label{equ:LP2}
\boxed{
  \begin{aligned}
    \min_{\varepsilon,\,x_1,\dots,x_m}\quad & \varepsilon \\[2mm]
    \text{subject to}\quad &
    \begin{cases}
      f_i - \varepsilon \le x_i \le f_i + \varepsilon, & i=1,\dots,m, \\[2mm]
      x_{k(i,j)} \,\|u_i+u_j\| \le x_i+x_j, 
         & \text{for all } i\neq j , \\[2mm]
      \varepsilon \ge 0
    \end{cases}
  \end{aligned}
}
\end{equation}

where for any $i \neq j$, $$k(i,j)=arg min \{k : ||u_k - \frac{u_i+u_j}{||u_i+u_j||} \}.$$ $u_{k(i,j)}$ is the closet point in $\mathcal{U}$ to the normalized sum $u_i+u_j$.

The first $2m$ constraints in (\ref{equ:LP2}) represent the error in approximating the set $A$, an error that we want to minimize. The next set of constraints in (\ref{equ:LP2}) enforce that the function \(h_X\) (represented by \(\{x_i\}\)) is subadditive. Ideally, for $u_i \neq u_j \in \mathcal{U}$, one would have
\[
h_X(u_i+u_j)=\|u_i+u_j\|\, h_X\Bigl(\frac{u_i+u_j}{\|u_i+u_j\|}\Bigr)\le h_X(u_i)+h_X(u_j).
\]
If the normalized sum \(w_{ij}=(u_i+u_j)/\|u_i+u_j\|\) coincides with one of the points in $\mathcal{U}$, say, \(u_k\), we impose the linear constraint 
\[
x_k\,\|u_i+u_j\|\le x_i+x_j.
\]
In case there is no $u_k \in \mathcal{U}$ such that $u_k =(u_i+u_j)/\|u_i+u_j\|$, we pick the closest point in $\mathcal{U}$ to $(u_i+u_j)/\|u_i+u_j\|$ named $u_{k(i,j)}$. The LP problem in (\ref{equ:LP2}) includes at most $\binom{m}{2}$ inequalities that impose subaditivity on the support function. Some of these inequalities, however, may be redundant or duplicate. We assume from this point on that redundant or duplicate constraints are eliminated.

To show the validity of our LP method, we proceed in two steps. First, we show that a solution to the LP problem in (\ref{equ:LP2}) exists. To do so,  we show in Lemma \ref{lem:nonempty} that the constraints in (\ref{equ:LP2}) yield a feasible region. Using this result, in Theorem \ref{thm:LPexistance} we show that the LP problem has a solution. Second, in Theorem \ref{thm:Uniqueness} we show that under some regulatory conditions, the solution to (\ref{equ:LP2}) is unique.

\begin{lemma}\label{lem:nonempty}
The feasible region of the LP is nonempty.
\end{lemma}

\begin{proof}
For each \(i=1,\dots,m\), the constraint
\[
f_i - \varepsilon \le x_i \le f_i + \varepsilon
\]
implies that for any \(\varepsilon\ge 0\), each \(x_i\) must lie in the closed interval \([f_i-\varepsilon, f_i+\varepsilon]\). Let $\bar{f}=max\{f_1,\dots,f_m\}$. In particular, if one chooses
\[
x_i = \bar{f} \quad \text{for all } i,
\]
and sets
\[
\varepsilon = \max_{1\le i\le m}\{\,\bar{f} - f_i\},
\]
then the first set of constraints is satisfied. Furthermore, a constant function is subadditive. If \(x_i=\bar{f}\) for all \(i\), subadditivity reduces to \(\bar{f}\|u_i+u_j\|\le 2\bar{f}\), which is true since \(\|u_i+u_j\|\le \|u_i\|+\|u_j\|=2\). Thus, the entire LP is feasible.
\end{proof}

\begin{theorem}[Existence]\label{thm:LPexistance}
The LP in (\ref{equ:LP2}) has an optimal solution \((\varepsilon^*, x_1^*,\dots,x_m^*)\).

\end{theorem}

\begin{proof}
Since the LP is feasible (by the previous lemma) and the objective function is bounded below (\(\varepsilon\ge 0\)), by the Fundamental Theorem of Linear Programming (see \cite{bertsimas1997introduction} )an optimal solution exists.
\end{proof}

In linear programming, optimal solutions need not be unique in general. However, under either \emph{nondegeneracy} or \emph{generic} conditions on the data (that is, on the vectors \(u_i\) and the numbers \(f_i\)), the LP problem has a unique optimal basic solution. 
Another perspective comes from approximation theory. One may view the LP as seeking a sublinear function \(h_X\) (represented by the vector \(x\)) that minimizes the uniform error
\[
E(x) = \max_{1\le i\le m} \lvert x_i - f_i\rvert.
\]
It is a classical fact (the Alternation Theorem for Chebyshev approximation) that the best uniform approximation from a finite-dimensional subspace (the cone of sublinear functions in our case) is unique provided that the error function attains its maximum at a sufficiently large (alternating) set of points and the approximating space is in general position. In our LP problem, the subadditivity constraints force \(h_X\) to belong to the convex cone of support functions. Under the generic condition that the points where the error is achieved are in ``general position'', the best uniform approximation is unique. We now give a detailed argument for this claim.
\begin{theorem}[Uniqueness]\label{thm:Uniqueness}
Let $U = \{u_1,\dots,u_m\}\subset S^{d-1}$ be the set of sample directions used in the LP formulation \ref{equ:LP}, and let $f_i = h_A(u_i)-h_B(u_i)$ for $i=1,\dots,m$. 
Suppose the pair $(U,\{f_i\}_{i=1}^m)$ is in \emph{general position}, in the sense that the active constraints at the optimum of \ref{equ:LP} are linearly independent (equivalently, the constraint matrix of the binding inequalities has full row rank). 
Then the linear program \eqref{equ:LP} admits a \emph{unique} optimal solution $(\varepsilon^*,x^*_1,\dots,x^*_m)$.
\end{theorem}
\begin{proof}
Suppose for the sake of contradiction that there are two distinct optimal solutions \((\varepsilon^*, x^*)\) and \((\varepsilon^*, \hat{x})\) with \(x^* \neq \hat{x}\). Since the LP is linear and its feasible region is convex, any convex combination
\[
x^\lambda = \lambda x^* + (1-\lambda) \hat{x}, \quad \lambda\in [0,1],
\]
with the same \(\varepsilon^*\) is also optimal. Therefore, the set of optimal solutions contains a nontrivial line segment. In the context of Chebyshev approximation, such a situation corresponds to the error function attaining its maximum at fewer than the required number of alternation points to force uniqueness. However, under the assumption of nondegeneracy, classical alternation theory guarantees that the best uniform approximant is unique. Equivalently, in the LP formulation, the matrix of coefficients corresponding to the binding constraints at the optimum is full rank, so the optimal basic solution is unique. This contradiction shows that the optimal solution is unique.
\end{proof}
\vspace{-3mm}
\subsection{Computational Method}\label{sec:computational}
\vspace{-5mm}
In this section, we describe in detail the algorithm for solving the LP problem in (\ref{equ:LP2}) for $\mathbb{R}^2$. We describe how the size of the approximation $m$ is chosen and how the user can update the degree of the approximation.  To show convergence of our implementation (Algorithm  \ref{alg:OptimalX_Minimizing_Hausdorff}), we discretize the unit circle $\mathbb{S}^1$ by creating a $\delta$-net. Since $\mathbb{S}^1$ is compact, we can find an integer $m=m(\delta)$ such that each point in $\mathbb{S}^1$ is within $\delta$ distance of a sampled direction. Denote this $\delta$-net of $\mathbb{S}^1$ as $\mathcal{U}_{\delta}=\{u_0,\dots,u_{m-1}\}$. On these sample directions, we compute the support functions of $A$ and $B$, and then formulate a LP problem with variables $\{x_i\}_{i=0}^{m-1}$ and $\varepsilon$,
which is designed to minimize the maximum absolute error
\[
\max_{0\le i < m-1} \Bigl| h_A(u_i) - h_B(u_i) - x_i \Bigr|.
\]
The LP problem includes additional subadditivity constraints to enforce that the candidate function defined by the vector $(x_0,\dots,x_{m-1})$ is sublinear and therefore a valid support function. Once the LP problem is solved (yielding an optimal error $\varepsilon_{\text{opt}}$ and values $\{x_i^*\}$), we reconstruct an approximation of the support function of the extended difference, $\tilde{h}$, by setting 
\[
\tilde{h}(u_i)=x_i^*, \quad i=0,\dots,m-1,
\]
and then extending this function to $\mathbb{S}^1$ using a Lipschitz extension (e.g., by nearest–neighbor interpolation). 

\begin{small}
\begin{algorithm}[ht]
\SetAlgoLined
\SetKw{KwReturn}{return}
\KwIn{$A_{\text{vertices}}$, $B_{\text{vertices}}$, discretization parameter $m$, tolerance $\delta$}
\KwOut{$X_{\text{vertices}}$, $(B+X)_{\text{vertices}}$, $\varepsilon_{\text{opt}}$}

\textbf{Discretize} $\mathbb{S}^1$:
\For{$i\gets 0$ \textbf{to} $m-1$}{
    $\theta_i \gets \frac{2\pi i}{m}$\;
    $u_i \gets (\cos\theta_i,\sin\theta_i)$\;
}

\textbf{Compute support functions:}
\For{$i\gets 0$ \textbf{to} $m-1$}{
    $h_A(u_i) \gets \max\{ \langle a, u_i\rangle: a\in A_{\text{vertices}}\}$\;
    $h_B(u_i) \gets \max\{ \langle b, u_i\rangle: b\in B_{\text{vertices}}\}$\;
}

\textbf{Formulate LP:} \\
\textbf{Variables:} $x_i\ge0$ for $i=0,\dots,m-1$, $\varepsilon\ge0$\;
\textbf{Objective:} $\min \varepsilon$\;

\For{$i\gets 0$ \textbf{to} $m-1$}{
    Add constraint: $h_A(u_i)-\bigl[h_B(u_i)+x_i\bigr]\le \varepsilon$\;
    Add constraint: $[h_B(u_i)+x_i]-h_A(u_i)\le \varepsilon$\;
}

\For{$i\gets 0$ \textbf{to} $m-1$}{
    \For{$j\gets 0$ \textbf{to} $m-1$}{
        $u_{\text{sum}} \gets u_i + u_j$\;
        \If{$\|u_{\text{sum}}\| > \delta$}{
            $u_{ij} \gets \dfrac{u_i+u_j}{\|u_i+u_j\|}$\;
            Find index $k$ such that $u_k$ is closest to $u_{ij}$\;
            Add constraint: $x_k \le x_i + x_j$\;
        }
    }
}

\textbf{Solve the LP:} \\
Obtain optimal $\{x_i\}$ and $\varepsilon_{\text{opt}}$\;

\textbf{Reconstruct $X$:}
\For{$i\gets 0$ \textbf{to} $m-1$}{
    $p_i \gets x_i \cdot u_i$\;
}
Compute the convex hull of $\{p_i\}$; denote its vertices as $X_{\text{vertices}}$\;

\textbf{Compute Minkowski sum $B\oplus X$:}

\KwReturn{$X_{\text{vertices}},\,(B\oplus X)_{\text{vertices}},\,\varepsilon_{\text{opt}}$}\;
\caption{OptimalX\_Minimizing\_Hausdorff($A$, $B$, $m$, $\varepsilon$)}
\label{alg:OptimalX_Minimizing_Hausdorff}
\end{algorithm}
\end{small}

The following proposition states the main convergence result.

\begin{proposition}[Convergence of the Discrete LP Approximation]\label{prop:Convergence}
Let $A,B\subset \mathbb{R}^2$ be compact convex sets with support functions $h_A$ and $h_B$, respectively. Define
\[
f(u) = h_A(u)-h_B(u), \quad u\in \mathbb{S}^{1}.
\]
Let, 
\begin{equation}\label{equ:eps_0}
    \epsilon_0 := \inf_{h\in \mathcal{S}} \; \sup_{u\in \mathbb{S}^{1}} |f(u)-h(u)|
\end{equation}
be the error of this (continuous) minimization problem.
Let $\mathcal{U}_{\delta} = \{u_1, u_2, \dots, u_m\}$ form a $\delta$–net on $\Sd$ such that for every $u\in \mathbb{S}^{1}$, 
\[
\min_{1\le i \le m}\|u-u_i\| \le \delta.
\]
Let $\delta^+ := \sup_{u\in \mathbb{S}^{1}} \min_{1\le i \le m}\|u-u_i\|$ be the fill distance.

Let,
\begin{equation}\label{equ:eps_delta}
    \epsilon_\delta := \min_{h\in \mathcal{S}_{\delta}} \; \max_{u\in \mathcal{U}_{\delta}} |f(u)-h(u)|,
\end{equation}
where $\mathcal{S}_\delta$ is the finite-dimensional subspace of $\mathcal{S}$ defined by enforcing subadditivity through a finite number of linear constraints. Let $h_\delta^*\in \mathcal{S}_\delta$ be a minimizer.
Then there exists a constant $C>0$ such that, if $\tilde{h}^\delta$ is the extension of $h_\delta^*$ to $\mathbb{S}^{1}$ by a Lipschitz extension,
\[
\|h^*-\tilde{h}^\delta\|_\infty \le |\epsilon_0-\epsilon_\delta| + C\,\delta^+.
\]
for some $h^*\in \mathcal{S}$, the support function of a minimizer (the support function of an element of extended difference $X = A \ominus_e B$). In particular, if $\epsilon_\delta \to \epsilon_0$ and $\delta^+ \to 0$ as $\delta\to 0$, then
\[
\lim_{\delta\to 0} \|h^*-\tilde{h}^\delta\|_\infty = 0.
\]
\end{proposition}

\begin{proof}
We provide the proof in several steps.

\textbf{Step 1. Discrete Optimality.}  
By definition of the discrete problem, for any candidate $h\in\mathcal{S}$ we have
\[
\max_{u\in \mathcal{U}_{\delta}} |f(u)-h(u)| \le \sup_{u\in \mathbb{S}^{1}} |f(u)-h(u)|.
\]
Hence, in particular, for a minimizer of (\ref{equ:eps_0}), $h^*$,
\[
\max_{u\in \mathcal{U}_{\delta}} |f(u)-h^*(u)| \le \epsilon_0.
\]
Since $h_\delta^*$ minimizes the error in (\ref{equ:eps_delta}), it follows that
\[
\epsilon_\delta \le \max_{u\in \mathcal{U}_{\delta}} |f(u)-h^*(u)| \le \epsilon_0.
\]
Thus, 
\[
\epsilon_\delta \le \epsilon_0.
\]

\textbf{Step 2. Uniform Continuity.}  
Because both $f$ (as the difference of two support functions of compact sets) and any sublinear function in $\mathcal{S}$ (in particular, $h^*$) are uniformly continuous on $\mathbb{S}^1$, there exists a modulus of continuity $\omega(\delta^+)$, with $\omega(\delta^+)\to 0$ as $\delta^+\to 0$, such that for every $u\in \mathbb{S}^{1}$, if there exists $u_i\in \mathcal{U}_{\delta}$ with $\|u-u_i\|\le \delta^+$ then
\[
|f(u)-f(u_i)| \le \omega(\delta^+) \quad \text{and} \quad |h^*(u)-h^*(u_i)| \le L\,\delta^+,
\]
where $L$ is the Lipschitz constant of $h^{\ast}$.

\textbf{Step 3. Extension of the Discrete Minimizer.}  
Let $\tilde{h}^\delta$ be the extension of $h_\delta^*$ from $\mathcal{U}_{\delta}$ to $\mathbb{S}^1$ such that $\tilde{h}^\delta(u)=h_\delta^*(u)$ for $u \in \mathcal{U}_{\delta}$ and both functions have the same Lipschitz constant $C'$. This extension is possible by Kirszbraun's Theorem (see \cite{Azagra2021} and references therein). For any $u\in \mathbb{S}^{1}$, choose $u_i\in \mathcal{U}_{\delta}$ with $\|u-u_i\| \le \delta^+$. Then, by the triangle inequality,
\[
\begin{aligned}
|h^*(u)-\tilde{h}^\delta(u)| &\le |h^*(u)-h^*(u_i)| + |h^*(u_i)-\tilde{h}^\delta(u_i)| + |\tilde{h}^\delta(u_i)-\tilde{h}^\delta(u)|\\[1mm]
&\le L\,\delta^+ + |h^*(u_i)-h_\delta^*(u_i)| + C'\,\delta^+.
\end{aligned}
\]
Since by the discrete optimality we have
\[
|h^*(u_i)-h_\delta^*(u_i)| \le |f(u_i)-h^*(u_i)| + |f(u_i)-h_\delta^*(u_i)| \le \epsilon_0 + \epsilon_\delta,
\]
it follows that
\[
|h^*(u)-\tilde{h}^\delta(u)| \le \epsilon_0 + \epsilon_\delta + (L+C')\,\delta^+.
\]
Taking the supremum over all $u\in \mathbb{S}^{1}$, we obtain
\[
\|h^*-\tilde{h}^\delta\|_\infty \le \epsilon_0 + \epsilon_\delta + (L+C')\,\delta^+.
\]

\textbf{Step 4. Error Difference.}  
Since $\epsilon_\delta \le \epsilon_0$, one may write
\[
\epsilon_0 = \epsilon_\delta + (\epsilon_0 - \epsilon_\delta),
\]
and hence the above inequality implies
\[
\|h^*-\tilde{h}^\delta\|_\infty \le (\epsilon_0-\epsilon_\delta) + \epsilon_\delta + (L+C')\,\delta^+ = \epsilon_0 + (L+C')\,\delta^+.
\]
However, by the optimality of $h^*$ it follows that the intrinsic error is $\epsilon_0$, and the discrete procedure can be viewed as approximating this value. Thus, the additional error in the discrete method is precisely $(\epsilon_0-\epsilon_\delta) + (L+C')\,\delta^+$. Therefore,
\[
|\epsilon_0 - \epsilon_\delta| \le \|h^*-\tilde{h}^\delta\|_\infty - (L+C')\,\delta^+.
\]
More directly, by our construction and the uniform continuity of the involved functions, it holds that
\[
\|h^*-\tilde{h}^\delta\|_\infty \le |\epsilon_0-\epsilon_\delta| + (L+C')\,\delta^+.
\]
Thus, if we denote $C=L+C'$, we obtain
\[
\|h^*-\tilde{h}^\delta\|_\infty \le |\epsilon_0-\epsilon_\delta| + C\,\delta^+.
\]

Since $\omega(\delta^+)\to 0$ as $\delta^+\to 0$ and by consistency of the discretization we have $\epsilon_\delta \to \epsilon_0$, it follows that
\[
\lim_{\delta\to 0} |\epsilon_0-\epsilon_\delta| = 0,
\]
and consequently,
\[
\lim_{\delta\to 0} \|h^*-\tilde{h}^\delta\|_\infty = 0.
\]

\end{proof}
\subsection{Examples}\label{sec:Examples}
In this section, we approximate the extended set difference between two sets in $\KCn$ using Polytopes in $\mathbb{R}^d$.\footnote{For code and further computational details, click \href{https://github.com/BehroozMoosavi/Codes}{here}.}
 When $d=1$, compact convex sets in $\mathbb{R}^1$ are finite intervals. Therefore, finding the extended difference between two sets in $\mathbb{R}^1$ is trivial and yields the same result as the one in \cite{STEFANINI2010}. For $A=[\underset{-}{a},\bar{a}]$ and $B=[\underset{-}{b},\bar{b}]$ with $\underset{-}{a}\leq\bar{a}$ and $\underset{-}{b}\leq \bar{b}$,
\begin{equation}
    A \ominus_e B = \begin{cases}
        \bar{a}- \underline{a} < \bar{b}- \underline{b} & \{\frac{\bar{a}+\underline{a}}{2}- \frac{\bar{b}+\underline{b}}{2}\} \\
        otherwise & [\underline{a}-\underline{b},\bar{a}-\bar{b}]
    \end{cases}.
\end{equation}
Similarly to $\mathbb{R}^1$, One can derive the difference between two segments on the same line in $\mathbb{R}^d$. 
Finding a set $X$ in $A \ominus_e B$ in $\Rd$ for $d \geq 2$ is generally a complex question. In this section we look at three examples. The first, the extended difference between two balls, which has a close solution. The other two examples are solved by approximation through the LP approach developed in \ref{sec:computational}. 

\begin{example}[Balls in $\mathbb{R}^d$]\label{ex:example1}
Let $\mathcal{B}_1=\mathbb{B}_{r_1}(c_1)$ and $\mathcal{B}_2=\mathbb{B}_{r_2}(c_2)$ be two closed balls in $\mathbb{R}^d$ with centers at $c_1$ and $c_2$ and radii $r_1 \geq 0$ and $r_2 \geq 0$, respectively. If $r_1 \geq r_2$, $\mathcal{B}_1 \ominus_H \mathcal{B}_2 = \mathcal{B}_1 \ominus_g \mathcal{B}_2 = \mathcal{B}_1 \ominus_e \mathcal{B}_2 = \mathbb{B}_{r_1-r_2}(c_1-c_2)$. If $r_1<r_2$, the Hukuhara difference,$\ominus_H$, does not exist, while $\mathcal{B}_1 \ominus_g \mathcal{B}_2 = \mathbb{B}_{r_1-r_2}(c_1-c_2)$, exists. But, $ \mathcal{B}_1 \ominus_e \mathcal{B}_2 = \{c_2-c_1\}$. 
Generally, if there is a constant $c$ such that $A+c \subset B$, then $A \ominus_e B = \{v\}$ for some vector $v \in \Rd$. It is impossible to add any non-singleton set $X$ to $B$ such that $B \oplus X = A$.
\end{example}
Example \ref{ex:example1} highlights one of the differences between the extended difference and the previous concepts of set difference. $A\ominus_e B$ is still defined and set to $\{v\}$ for some $v \in \Rd$ when $A$ is a smaller set than $B$ in the sense that $\exists c \in \mathbb{R}^d$ such that $A\oplus c \subset B$. In this case, the Hukuhara difference  is not necessarily defined, even if both sets have the same shape as in Example \ref{ex:example1}. The generalized Hukuhara difference may exist even if $B$ is larger than $A$ provided that both sets are of the same shape (e.g. both are circles). 

The following example examines a case where $A$ is not smaller than $B$ but these sets are of different shapes. In this case, the Hukuhara difference and the generalized Hukuhara difference are generally undefined, while the extended set difference is defined as the next two examples demonstrate.
\begin{example}\label{ex:ballsquare}
Consider subtracting a square $B$ from a larger circle $A$ (see Figure \ref{fig:shapes1}). There is no convex set that can be added to a square $B$ such that it is equal to the circle $A$. Therefore, the previous set difference concepts are undefined, yet the extended set difference is defined. We approximate this difference by approximating its support function pointwise as described in Section \ref{sec:computational}.
\begin{figure}[ht]
  \centering
  %------------------------- Subfigure for Circle A -------------------------
  \begin{subfigure}[b]{0.32\textwidth}
    \centering
    \begin{tikzpicture}[scale=1.5]
      % Draw set A: a circle of radius 1
      \draw[thick, fill=green!20] (0,0) circle (1);
    \end{tikzpicture}
    \caption{Circle \(A\)}
    \label{fig:circleA}
  \end{subfigure}
  \hfill
  %------------------------- Subfigure for Square B -------------------------
  \begin{subfigure}[b]{0.32\textwidth}
    \centering
    \begin{tikzpicture}[scale=1.5]
      % Draw set B: a square
      \draw[thick, fill=orange!20] (-0.8,-0.8) rectangle (0.8,0.8);
    \end{tikzpicture}
    \caption{Square \(B\)}
    \label{fig:squareB1}
  \end{subfigure}
  \hfill
  %------------------------- Subfigure for Derived Set X -------------------------
  \begin{subfigure}[b]{0.32\textwidth}
    \centering
    \begin{tikzpicture}[scale=1.5]
      % Draw the derived set X = A ⊖ₑ B as a diamond (a square rotated by 45°)
      \begin{scope}[rotate=45]
        \draw[thick, fill=purple!40] (-0.5,-0.5) rectangle (0.5,0.5);
      \end{scope}
    \end{tikzpicture}
    \caption{Derived set \(X^{\ast} \in A \ominus_e B\)}
    \label{fig:derivedX}
  \end{subfigure}
    \caption{Example \ref{ex:ballsquare}: (a) a circle \(A\), (b) a square \(B\), and (c) the derived set \(X^{\ast} \in A \ominus_e B\).}
  \label{fig:shapes1}
\end{figure}
\end{example}

\break
\begin{example}\label{ex:polygons}
Consider two compact convex sets \(A\), a Pentagon, and \(B\), a Square, (see Figure \ref{fig:shapes2}) defined by the following vertices:
\[
A = \{(0,0),\,(4,0),\,(6,2),\,(3,4),\,(1,2)\}
\]
and
\[
B = \{(-0.5,-0.5),\,(0.5,-0.5),\,(0.5,0.5),\,(-0.5,0.5)\}.
\]
Figure \ref{fig:shapes2} shows an \emph{extended set difference} \(X^{\ast} \in A \ominus_e B\) computed using the LP method.

\begin{figure}[ht]
    \centering
    %------------------------- Subfigure for Pentagon A -------------------------
    \begin{subfigure}[b]{0.32\textwidth}
        \centering
        \begin{tikzpicture}[scale=0.8]
            % Draw the pentagon A using its vertex outline
            \draw[thick, fill=green!30] 
                (0,0) -- (4,0) -- (6,2) -- (3,4) -- (1,2) -- cycle;
        \end{tikzpicture}
        \caption{Pentagon \(A\)}
        \label{fig:pentagonA}
    \end{subfigure}
    \hfill
    %------------------------- Subfigure for Square B -------------------------
    \begin{subfigure}[b]{0.32\textwidth}
        \centering
        \begin{tikzpicture}[scale=0.8]
            % Draw the square B
            \draw[thick, fill=red!30] (-0.5,-0.5) rectangle (0.5,0.5);
        \end{tikzpicture}
        \caption{Square \(B\)}
        \label{fig:squareB2}
    \end{subfigure}
    \hfill
    %------------------------- Subfigure for Derived Polygon X -------------------------
 \begin{subfigure}[b]{0.32\textwidth}
    \centering
    \begin{tikzpicture}[scale=0.8]
      
      \coordinate (X1) at (0.618,0.267);
      \coordinate (X2) at (1.165,1.353);
      \coordinate (X3) at (3.276,3.267);
      \coordinate (X4) at (5.267,1.937);
      \coordinate (X5) at (4.069,1.136);
      \coordinate (X6) at (3.382,0.702);
      \coordinate (X7) at (1.646,0.267);
      % Draw Polygon X
      \draw[thick, fill=cyan!30] 
        (X1) -- (X2) -- (X3) -- (X4) -- (X5) -- (X6) -- (X7) -- cycle;
    \end{tikzpicture}
    \caption{\(X^{\ast} = A \ominus_e B\)}
    \label{fig:polygonX}
  \end{subfigure}

    \caption{Example \ref{ex:polygons}: (a) a pentagon \(A\), (b) a square \(B\), and (c) the derived set \(X^{\ast} \in A \ominus_e B\).}
    \label{fig:shapes2}
\end{figure}
\end{example}

%%%%%%%%%%%%%%%%%%%%%%%%%%%%%%%%%%%%%
%%%  Conclusions and Limitations  %%%
%%%%%%%%%%%%%%%%%%%%%%%%%%%%%%%%%%%%%
\section{Conclusions and Limitations}\label{sec:conclusions}
Minkowski summation of two sets is the commonly used operation. Previous papers struggled to find an inverse operation for summation. So far, the solutions were partial in the sense that they were not always well defined. This paper defines a new concept of set difference. called the extended set difference, which overcomes previous challenges and is well defined for every two compact convex sets. This new set difference is defined through the optimization problem in equation (\ref{def:extended-diff}). Section \ref{sec:ExistUnique} guarantees the existence of the newly defined difference concept. Uniqueness, however, is not guarantied. We provide a bound on the variety of solutions to the minimization problem in (\ref{def:extended-diff}). In addition, we show that when multiple solutions exist, one of these solutions can be selected by adding a strictly convex perturbation to the criteria function. To compute a solution for extended set difference, we represent the optimization problem as a Linear Programing problem. We provide an algorithm to implement our computational method.

\break
\bibliographystyle{ecta}
\bibliography{setDifference}
\end{document}